\documentclass[12pt]{article} 

\usepackage{amsmath}
\usepackage[dvips]{graphicx}
\usepackage{amsfonts}
\usepackage{amsopn}
\usepackage{amsthm}

\newcommand{\Ob}{\bar{\Omega}}
\renewcommand{\O}{\Omega}
\newcommand{\Tb}{\bar{T}}
\newcommand{\T}{T}

\newcommand{\e}{\epsilon}

\renewcommand{\a}{\alpha}
\renewcommand{\b}{\beta}
\newcommand{\z}{\zeta}
\renewcommand{\o}{\omega}
\newcommand{\Pt}{\tilde{P}}
\newcommand{\F}{\mathcal{F}}

\renewcommand{\epsilon}{\varepsilon}

\renewcommand{\)}{\right)}

\newcommand{\BI}{B_{1}}
\newcommand{\BII}{B_{2}}
\newcommand{\BIII}{B_{3}}

\renewcommand{\l}{\mathcal{L}}

\DeclareMathOperator{\supp}{supp}
\DeclareMathOperator{\Lip}{Lip}
\DeclareMathOperator{\diam}{diam}

\newtheorem{theorem}{\bf \normalsize \bf  Theorem}[section]

\newtheorem{corollary}[theorem]{\bf \normalsize \bf Corollary}
\newtheorem{lemma}[theorem]{\bf \normalsize \bf Lemma}
\newtheorem{proposition}[theorem]{\bf \normalsize \bf Proposition}

\newtheorem{problem}{\bf \normalsize \bf  Problem}
\newtheorem*{problem'}{\bf \normalsize \bf  Problem I$_{\text{\bf w}}$}

\theoremstyle{remark}
\newtheorem{remark}[theorem]{\bf \normalsize \bf Remark}

\theoremstyle{definition}
\newtheorem{definition}[theorem]{\bf \normalsize \bf Definition}

\date{}
\title{Optical design of two-reflector systems,  the
Monge-Kantorovich mass transfer  problem and Fermat's principle}
\author{Tilmann Glimm \& Vladimir Oliker
\\
\\
Department of Mathematics and Computer Science,\\
 Emory University, Atlanta, Georgia 30322\\
}
\date{May 16, 2002}
\begin{document}

\maketitle
\begin{abstract}
It is shown that the problem of designing a two-reflector system
transforming a plane wave front with given intensity into
an output plane front with prescribed output intensity can be 
formulated and solved as the Monge-Kantorovich  mass transfer problem.  

\end{abstract}

\section{Introduction} 

\begin{figure}
\begin{center}
\includegraphics[width=12cm]{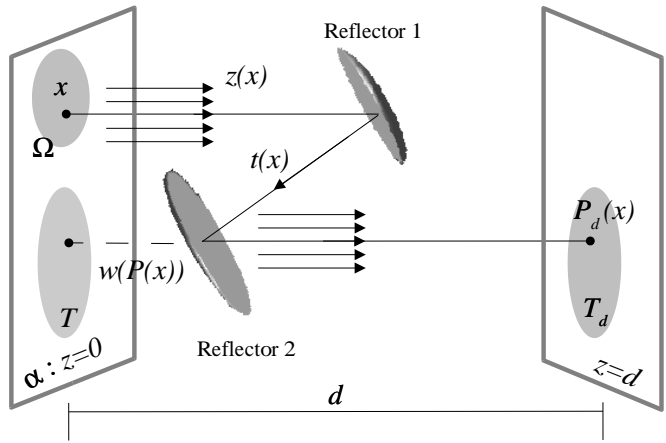}
\caption{Sketch for Problem \ref{pmain}} \label{figure geometr opt}
\end{center}
\end{figure}

Consider a two-reflector system of configuration shown schematically on 
Fig. \ref{figure geometr opt}. Let $(x = (x_1,x_2,...,x_n),z)$ be
the Cartesian coordinates in ${\mathbb R}^{n+1},~n \geq 2,$ with 
$z$ being the horizontal axis and $x_1,x_2,...,x_n$ the coordinates
in
the hyperplane $\a: z=0$. Let $\BI$ denote a beam of parallel 
light rays propagating in the positive
$z-$direction and let $\Ob$ denote the wavefront which is
the cross section of $\BI$ by hyperplane $\a$.
Assume that $\O$ is a bounded domain on $\a$.
An individual ray of the front is labeled by a point $x \in \Ob$. 
The light intensity of the beam $\BI$ is
denoted by $I(x), x \in \O$, where $I$ is a non-negative integrable function.

The incoming beam  $\BI$ is intercepted by  the first reflector 
$R_1$, defined as a graph of a function $z(x), x \in \Ob$. The rays in $\BI$
are reflected off $R_1$ 
forming a beam of rays $\BII$. The beam $\BII$ is intercepted by
reflector $R_2$ which transforms it into the output beam $\BIII$. The 
beam   $\BIII$ also consists of parallel light rays propagating in the 
same direction as $\BI$. The output wavefront at a distance $d >0$ from
the hyperplane $\a$ is denoted by  $\Tb_d$; the projection of $\Tb_d$ on the
hyperplane $\a$  we denote by  $\Tb$. The second reflector $R_2$ is also assumed 
to be a
graph of a function $w(p), p\in\Tb$. The quantity
$\frac{1}{|J(P_d(x)|},$
where $P_d$ is the map of $\Ob$ on $\Tb_d$ and $J$ is the Jacobian, 
is the {\it {expansion ratio}} 
and it measures the expansion of a tube of rays 
due to the two reflections \cite{Keller:95}. 
It is assumed that both $R_1$ and $R_2$ are perfect
reflectors and no energy is lost in the transformation process.
Consequently,  the corresponding relation between the 
input intensity $I$ on $\O$ and output intensity $L$ on $T_d$
is given by 
\begin{equation} \label{ec0}
L(P_d(x))|J(P_d(x)| = I(x).
\end{equation}

The ``two-reflector'' problem that needs to be solved by designers of optical systems
consists in determining the reflectors $R_1$ and $R_2$ so that all of the
properties of the two-reflector system above hold 
for prescribed in advance domains 
$\O, T$ and
positive integrable functions $I(x), x \in \O,$ and $L(p), p \in T$; 
see Malyak \cite{Malyak} and other references there. 
It is usually assumed 
in applications that $\O$ and $T_d$ are bounded and convex.

Two fundamental principles of geometrical optics are used
to describe the transformation of 
the beam  $\BI$ into beam $\BIII$: 
the classical reflection law leading to the ray tracing equations
defining the map $P_d$ 
and the energy conservation law for the energy flux along infinitesimally
small tubes of 
rays; see \cite{Malyak} where the problem is formulated for 
rotationally symmetric data and a class of rotationally symmetric 
solutions is found. 

The problem of recovering reflectors $R_1$ and $R_2$ without assuming
rotational symmetry
was formulated rigorously by Oliker and Prussner in \cite{ACES},
and it was shown that it can be considered as a problem of determining
a special map of $\Ob \rightarrow \Tb$ with a potential 
satisfying an equation of 
Monge-Amp\`ere type relating the input and output intensities. 
Existence and
uniqueness of weak solutions were established by Oliker at that time
but only the numerical results implementing a constructive 
scheme for proving existence
were presented in \cite{ACES} for several test cases. Detailed proofs were
given in \cite{Oliker:TiNA}. 

In this paper we show that this problem can also be studied in the
framework of the Monge-Kantorovich mass transfer problem studied by Brenier \cite{Brenier1},   
Caffarelli \cite{Caf_alloc:96}, Gangbo and McCann \cite{Gangbo/McCann:95}, and
other authors. In our notation, the  Monge-Kantorovich mass transfer problem 
is
to transfer the intensity  $I$ on $\O$ into the intensity $L$ on $T$ 
via a map $P:\O \rightarrow \T$ for which
the total transportation cost $\int_{\O} C(x,P(x))\,Idx$ is minimized.
Here $C(x,p)$ is a given strictly convex cost function.

\begin{sloppypar}
The proof of existence and uniqueness of solutions
to  the Monge-Kantorovich problem is obtained by solving
a minimization problem for the functional
\end{sloppypar}
\begin{equation} \label{intro eq 2}
(\zeta,\o)\mapsto \int_{\O}\zeta Idx-\int_{T}\o Ldp
\end{equation}
considered on pairs of continuous functions 
$\zeta$ on $\Ob$ and $\o$ on $\Tb$ that satisfy
\begin{equation} \label{intro eq 1}
\zeta(x)-\o(p)\geq -C(x,p),~~x \in \Ob, ~p \in \Tb.
\end{equation}
Under various conditions it is shown 
in \cite{Brenier1}, \cite{Caf_alloc:96}, \cite{Gangbo/McCann:95}
that this functional
is minimized by some pair $(\zeta_0,\o_0)$ (referred to as 
Kantorovich potentials), and that 
$P(x)=x+\nabla \zeta_0$ solves the Monge-Kantorovich problem.

Applying these ideas, we show that the geometric optics problem at hand
can be formulated as a Monge-Kantorovich mass transfer problem
with a quadratic cost function; see section \ref{MK}.
The Kantorovich
potentials correspond to the pair of reflectors that solve the problem. 
The condition (\ref{intro eq 1})
has a geometric meaning; namely, it filters out 
 reflectors that allow only 
optical paths longer than a certain prescribed one. The functional 
(\ref{intro eq 2}) to be minimized
is the mean horizontal distance between points of the two reflectors, 
with the average weighted by the
two intensities.

We prove that there are always two different reflectors system satisfying
the stated requirements.   The corresponding ways in which 
one intensity is transfered into the
other one are exactly the {\it most} and the {\it least} energy efficient   
in the sense of the Monge-Kantorovich 
cost. This result is thus ultimately a 
variant of Fermat's principle.

The fact that the solution to the above geometrical optics problem
 can be derived 
from a variational
principle gives rise to a numerical treatment of the problem different 
from the one suggested in
\cite{Oliker:TiNA}. In particular, when the problem of minimizing
(\ref{intro eq 2}) under constraints (\ref{intro eq 1}) is discretized
we have the linear programming problem with quadratic constraints.
We ran some numerical experiments with this approach 
and intend to return to this point in a separate publication.

This paper is organized as follows. 
In section \ref{stpr}, we recall
some results from \cite{Oliker:TiNA} concerning the ray tracing map, 
assuming smoothness of the reflectors, and formulate the main 
``two-reflector''
problem.
In section \ref{geom}, 
we give a geometric characterization of reflectors as envelopes of certain
families of paraboloids. Such  characterization is of independent interest.
 In section \ref{weak problem} 
we use this geometric characterization to define weak solutions of type A
and type B of
the two-reflector problem. To prove existence and uniqueness of solutions
for each type we utilize the ideas of the 
Monge-Kantorovich theory and  introduce the functional (\ref{intro eq 2})
on a certain class of ``quasi-reflector'' systems. This is done in section
\ref{min_functional}. In the same section it is shown that the problem
of finding weak solutions of type A is equivalent to finding
minimizers of (\ref{intro eq 2}). Weak solutions of type B correspond
to maximizers of (\ref{intro eq 2}). On the other hand, existence of
minimizers (maximizers) to this functional is not difficult and has
been established before in
 \cite{Brenier1}, \cite{Caf_alloc:96}, \cite{Gangbo/McCann:95}.
This implies existence of solutions. The uniqueness in respective class
is established in section  \ref{MK} 
by proving that the ray tracing map $\Pt$ associated with a weak solution
minimizes or maximizes the quadratic Monge-Kantorovich cost for
which the functional (\ref{intro eq 2}) is the dual. The main theorem
on existence and uniqueness of weak solutions to the two-reflector
problem is stated and proved in section \ref{proof}.

Finally, we note that same methods lead to variational
formulations and solutions of other geometrical optics problems involving 
systems with single and multiple reflectors.

\section{Statement of the problem} \label{stpr}

We begin by reviewing briefly the analytic formulation of
the problem for smooth reflectors;  see \cite{Oliker:TiNA} for more details.

Let $R_1$ be given by the position vector
$r_1(x) = (x,z(x)),~ x \in \Ob,$ with $z \in C^2(\Ob)$.
The unit normal $u$  on
$R_1$ is given by 
\[
u = \frac{(-\nabla z,1)}{\sqrt{1 +
|\nabla z|^2}}.
\]
Consider a ray labeled by
$x \in \Ob$ and propagating in the positive direction $k$ of the $z-$ axis. 
According to the reflection law
the direction of the ray $y(x)$ reflected off $R_1$ is given by 
\[y = k - 2\langle k,u \rangle u = 
k - 2\frac{(-\nabla z,1)}{1 +|\nabla z|^2},\]
where $\langle, \rangle$ is the inner product in $R^{n+1}$.
Denote by $t(x)$ the distance from
reflector $R_1$ to reflector $R_2$ along the ray reflected
in the direction
$y(x)$ and let $s(x)$ be the distance from $R_2$ to the wavefront $\Tb_d$ along the
corresponding ray reflected off $R_2$.  Assume for now
that $t \in C^1(\Ob)$ and $R_2$ is a $C^1$ hypersurface.
The total optical path length
(OPL) corresponding to the ray associated with the point $x \in \Ob$ is
$l(x)= z(x) + t(x) + s(x).$
A calculation shows that $l(x)  = const \equiv l~\mbox{on}~ \Ob$.
Since
\begin{equation}
R_2:~~~r_2(x) = r_1(x) + t(x)y(x), ~~x \in \Ob,
\label{R_2}
\end{equation}
the image of $x$ on the reflected wavefront $\Tb_d$ is given by
\begin{equation}
P_d(x)=r_1(x) + t(x)y(x) +s(x)k,~~~x \in \Ob.
\label{P}
\end{equation}
The equation (\ref{P}) is the ray tracing 
equation for this two-reflector system. 

Introduce the map $P(x) = P_d(x) -dk: \Ob \rightarrow \Tb$. A
calculation \cite{Oliker:TiNA} shows  that
\begin{equation} \label{splitx1}
p = P(x) = x + \b \nabla z(x),~~x \in \Ob,
\end{equation}
where $\b = l - d$ is the ``reduced'' optical path length.

To simplify the notation we will write $L(P(x))$ instead of 
$L(P_d(x)) (\equiv L(P(x)+dk))$.
For the input intensity $I(x),~ x \in \O,$ and the
output intensity $L(P(x))$ on $T_d$ we have in 
accordance with the differential form of the energy conservation law (\ref{ec0}),
\begin{equation} \label{econ}
L(P(x))|J(P(x))| = I(x),~~x \in \O,
\end{equation}
where we also take into account that $J(P_d) =J(P)$. It follows from (\ref{econ}) that $\O,~T,~I$  and $L$
must satisfy the necessary condition
\begin{equation} \label{EC}
\int_T L(p)dp = \int_\O I(x)dx.
\end{equation}

It follows from (\ref{splitx1}) that 
$J(P)=\mbox{det}\left[Id + \b Hess(z)\right]$, 
where $Id$ is the identity matrix and $Hess$ is the Hessian. 
Hence, by (\ref{econ}),
\begin{equation} \label{ME}
L(x + \b \nabla z)|\mbox{det}\left[Id + \b Hess(z) \right]|
= I  ~~\mbox{in}~~ \O.
\end{equation}
Thus, the problem of determining the reflectors $R_1$ and $R_2$ 
with properties described in the introduction requires
solving the following
\begin{problem} \label{pmain}
Given bounded domains $\O$ and $\T$ on the hyperplane $\a$ 
and two nonnegative, 
integrable functions $I$ on $\O$
and $L$ on $\T$ satisfying (\ref{EC}), it is required to find a function 
$z \in C^2(\Ob)$ such that the
map
\begin{equation}\label{bc}
P_{\a}= x + \b \nabla z: \Ob \rightarrow \Tb
\end{equation}
is a diffeomorphism satisfying equation (\ref{ME}). 
\end{problem}
It is shown in \cite{ACES} that once such a function $z$ is found,
the function $w$ describing the
second reflector is determined
by $z$ and $\beta$ as
\begin{equation} \label{eqw}
w(P(x))=d-s(x)=z(x)+\frac{\beta}{2}(|\nabla z|^2-1).
\end{equation}

Following \cite{Oliker:TiNA} we introduce the function  
\begin{equation}\label{potential}
V(x) = \frac{x^2}{2} + \b z(x)-\frac{\b^2}{2}.
\end{equation}
Then by (\ref{splitx1}) and (\ref{eqw}) 
\begin{eqnarray}
P = \nabla V, \label{legendre1}\\
w= \frac{1}{\b}\left [ V - \langle x, \nabla V \rangle + 
\frac{1}{2}|\nabla V|^2 \right ], \label{legendre2}
\end{eqnarray}
where $V - \langle x, \nabla V \rangle$ is the negative of the usual 
Legendre 
transform of $V$. Thus, $V$ is a potential for
the map $P: \Ob \rightarrow \Tb$. If $P$ is
a diffeomorphism then the inverse of the transformation 
$(x,z(x)) \rightarrow (p, w(p))$, where $p = P(x)$, is given by
\begin{eqnarray} \label{invP}
x(p) = P^{-1}(p)=p-\beta \nabla_p w(p), \label{invleg1}\\
z(p) = w(p) - \frac{\b}{2}|\nabla_pw(p)|^2
+ \frac{\b}{2}),~~~p \in \Tb. \label{invleg2}
\end{eqnarray}

In terms of the potential $V$ the equation (\ref{ME}) becomes
\begin{equation} \label{ME1}
L( \nabla V)|\mbox{det}Hess(V)|= I  ~~\mbox{in}~~ \O,
\end{equation}
which is an equation of Monge-Amp\`{e}re type. 


In order to clarify the relations between the parameters $l$, $d$ and 
$\beta$, note first that
it is the reduced optical path
length $\beta$ that is intrinsic to the problem. 
The choice of  cross sections of the fronts 
(that is, the  selection of a particular value for 
$d$) is extraneous. In fact, it is easy to see that if $(z,w)$ are two reflectors
as above and 
we change $d$ to $d'$ then $l'-d' = l-d =\b$ and $(z,w)$ are not
affected by such change.

Finally, we note that the two-reflector system described above 
has the following two symmetries.
For $\lambda \in {\mathbb R}^+$ put $z'(x)= \lambda z(x)$ and 
$\beta'=\frac{1}{\lambda}\beta$. Then 
\begin{eqnarray*}
P'(x) = P(x), \\
w'(p)= \lambda w(p)+\frac{\beta}{2}(\lambda-\frac{1}{\lambda}),
\end{eqnarray*}
where $P'(x) = x + \b ' \nabla z'(x)$. In other words,
the system is invariant under some combination of flattening (stretching) 
the first and translating and flattening (stretching) the second reflector.

Note also that a horizontal translation of both reflectors, that is,
 adding the same
constant to $z$ and $w$, does not change $\b$ and the map $P$.

\section{Geometric characterization of reflectors} \label{geom}

We examine first more closely the relationship between the functions $z, V$
and $w$ for smooth reflectors. 
{\it {For the rest of the paper
we assume that $\O$ and $T$ are bounded domains on the hyperplane $\a$.}} 
We continue to assume that
the map $P: \Ob \rightarrow \Tb$ is a diffeomorphism. 
Let $x \in \Ob$, $p \in \Tb$
and 
\[Q(x,p) = \langle x, p \rangle + \b w(p) - \frac{p^2}{2}.\]
If $p = P(x)$ then by (\ref{potential}) - (\ref{legendre2}) 
we have
\[V(x) = Q(x,P(x)).\]
Denote by $S_V$ the graph of $V$ over $\Ob$.
Let $x_0 \in \Ob$ and $p_0 = P(x_0)$. The tangent hyperplane to $S_V$ 
at $(x_0, V(x_0))$ is given by the equation 
\[Z= \langle x, p_0 \rangle - \langle x_0, p_0 \rangle + V(x_0),\]
where $(x,Z)$ denotes  an arbitrary point on that hyperplane. 
Taking into account (\ref{legendre2}), we obtain
\begin{equation} 
\langle x, p_0 \rangle - \langle x_0, p_0 \rangle + V(x_0)= 
 \langle x, p_0 \rangle + \b w(p_0) - \frac{p_0^2}{2} = 
Q(x,p_0). 
\end{equation}
Since $V(x_0)= Q(x_0,p_0)$, we conclude that $Z = Q(x,p_0)$ is the
tangent hyperplane to $S_V$ at $(x_0, V(x_0))$.
Consequently, if $V$ is convex then $Z = Q(x,p_0)$ is a supporting
hyperplane to $S_V$ from below and if $V$ is concave then 
$Z = Q(x,p_0)$ is supporting to $S_V$ from above (relative to positive direction
of the $z$-axis). 
Because $\Tb$ is bounded there are no vertical tangent
hyperplanes to the graph of $V$ and we have
\begin{eqnarray}
V(x) \geq Q(x,p_0)~\mbox{for all}~x \in \Ob~~
\mbox{if $V$ is convex}, \nonumber \\
V(x) \leq Q(x,p_0) ~\mbox{for all}~x \in \Ob~~
\mbox{if $V$ is concave} \nonumber.
\end{eqnarray}
Since for every $p \in \Tb$ the hyperplane $Q(x',p)$ is 
supporting to $S_V$ at some $(x',V(x'))$, we get
\begin{eqnarray}
V(x) \geq Q(x,p)~\mbox{for all}~x \in \Ob,~p \in \Tb~~
\mbox{if $V$ is convex}, \label{convexV}\\
V(x) \leq Q(x,p) ~\mbox{for all}~x \in \Ob,~p \in \Tb~~
\mbox{if $V$ is concave} \label{concaveV}
\end{eqnarray}
and in both cases we have equalities if $p=P(x)$.

Let 
\begin{equation}
U(p)= \frac{p^2}{2} - \b w(p) - \frac{\b ^2}{2},  ~
R(x,p) = \langle x, p \rangle - \b z(x) - \frac{x^2}{2}. \nonumber
\end{equation}
It follows from (\ref{convexV}) and (\ref{concaveV}) that
\begin{eqnarray}
U(p) \geq R(x,p)~\mbox{for all}~x \in \Ob,~p \in \Tb~~
\mbox{if $V$ is convex}, \label{convexU}\\
U(p) \leq R(x,p) ~\mbox{for all}~x \in \Ob,~p \in \Tb~~
\mbox{if $V$ is concave} \label{concaveU}
\end{eqnarray}
and in both cases  equalities are achieved if $p=P(x)$. Also,
for any fixed $x_0 \in \Ob$ the hyperplane $ R(x_0,p)$ is
supporting to the graph $S_U$ of $U$ at $(p_0=P(x_0), U(p_0))$.

Using the usual characterization of convex functions 
\cite{Schneider} we obtain from (\ref{convexV}),  (\ref{convexU}) and
(\ref{concaveV}),  (\ref{concaveU})
\[V(x) = \sup_{p \in \Tb}Q(x,p), ~~
U(p) = \sup_{x \in \Ob}R(x,p)~~ \mbox{when $V$ is convex},\]
\[V(x) = \inf_{p \in \Tb}Q(x,p), ~~
U(p) = \inf_{x \in \Ob}R(x,p)~~ \mbox{when $V$ is concave}.\]
For convex $V$ this implies
\begin{eqnarray}
z(x)=\sup_{p\in\Tb}\left [ \frac{\beta^2-|x-p|^2}{2\beta}
+w(p)\right ],\enspace x\in\Ob.\label{convz}\\
w(p) = \inf_{x\in\Ob}\left[ \frac{|x-p|^2-\beta^2}{2\beta}
		+z(x)\right ],
\enspace p\in\Tb.\label{convw}
\end{eqnarray}
Similarly, when $V$ is concave we have
\begin{eqnarray}
z(x)=\inf_{p\in\Tb}\left [ \frac{\beta^2-|x-p|^2}{2\beta}
+w(p)\right ],\enspace x\in\Ob.\label{concz}\\
w(p) = \sup_{x\in\Ob}\left[ \frac{|x-p|^2-\beta^2}{2\beta}
		+z(x)\right ],
\enspace p\in\Tb.\label{concw}
\end{eqnarray}

The characterizations  (\ref{convz}), (\ref{convw}) and 
(\ref{concz}), (\ref{concw}) have a 
simple geometric meaning. To describe it, consider first
the case when $V$ is convex. Recall that the
total optical path length $l = z(x) + t(x) + d - w(P(x)) = const$ (see
Fig. 1). Also,
$t^2(x) = |x-P(x)|^2 + |z(x)-w(P(x)|^2$ and $\b = l-d$. It follows from
(\ref{convz}) 
that $(x,z(x))$ is a point on the graph of the paraboloid
\begin{equation} \label{eqk}
k_{p,w}(x)=\frac{\beta^2-|x-p|^2}{2\beta}
		+w,\enspace x\in\a,
\end{equation}
with the focus at $(p=P(x),w(P(x)))$ and focal parameter $\b$.

Similarly, it follows from (\ref{convw}) that 
a point $(p=P(x),w(P(x)))$ on the second reflector lies
on a paraboloid
\begin{equation}  \label{eqh}
h_{x,z}(p)=\frac{|x-p|^2-\beta^2}{2\beta}
		+z, \enspace p\in\a,
\end{equation}
with the focus at $(x,z(x))$ and focal parameter $\b$.

Let $K_{p,w(p)}$ be the convex body bounded by the graph of paraboloid
$k_{p,w}(x)$ and $H_{x,z(x)}$ the convex body bounded by the graph of paraboloid
$h_{x,z}(p)$. Then (\ref{convz}) and (\ref{convw}) mean that  the graphs $S_z$ 
of $z(x)$ and $S_w$  of $w(p)$ are given by
\begin{eqnarray}
        S_z =\partial \left (\bigcup_{p\in \Tb} K_{p,w(p)}\right),\label{convz1} \\
       S_w =\partial \left(\bigcup_{x\in \Ob} H_{x,z(x)}\right). \label{convw1}
\end{eqnarray}

When the potential $V$ is concave we have similar characterizations
of  $S_z$ and $S_w$ with $\bigcup$ in (\ref{convz1}),
(\ref{convw1}) replaced by $\bigcap$.

\begin {remark} It follows from (\ref{convz}), (\ref{convw})
 that when 
$V(x)$ is convex then for any $x \in \O$ 
the path taken by the light ray through
the reflector system is the shortest among all possible 
paths (not necessarily satisfying the reflection law)
 that go from $(x,0)$ to $(x,z(x))$, then to
some $(p,w(p))$ and then to $(p,d)$. Of course, the shortest 
path satisfies the reflection law and $p = P(x)$.
For concave $V$ the corresponding light path is the longest as it
follows from (\ref{concz}), (\ref{concw}). Thus the characterizations
(\ref{convz}), (\ref{convw}) and (\ref{concz}), (\ref{concw}) are
variants of the 
Fermat principle. 
\end{remark}

\section{Weak solutions of Problem \ref{pmain}} \label{weak problem}

We use the geometric characterizations of reflectors in
section \ref{geom} to define weak solutions to
Problem \ref{pmain}.  
Let $\O$ and $T$ be two bounded domains on the hyperplane $\a$ 
and $\b$ a fixed positive number.

\begin{definition} \label{type A-B}
A pair $(z,w)\in C(\Ob)\times C(\Tb)$ is called a  
{\em two-reflector of  type A} if
\begin{eqnarray} \label{type A}
z(x)=\sup_{p\in\Tb}k_{p,w(p)}(x),~ x \in \Ob \label{A1},\\
w(p)=\inf_{x\in\Ob}h_{x,z(x)}(p), ~p \in \Tb, \label{A2}
\end{eqnarray}
where $k_{p,w(p)}(x)$ and $h_{x,z(x)}(p)$ 
are defined by (\ref{eqk}) and (\ref{eqh}). 
Similarly, a pair $(z,w)\in C(\Ob)\times C(\Tb)$ is called a  
{\em two-reflector of  type B} if
\begin{eqnarray} \label{type B}
z(x)=\inf_{p\in\Tb}k_{p,w(p)}(x),~ x \in \Ob \label{B1},\\
w(p)=\sup_{x\in\Ob}h_{x,z(x)}(p), ~p \in \Tb.\label{B2}
\end{eqnarray}
\end{definition}

To avoid repetitions, we consider below only two-reflectors of type
A. The changes that need to be made to deal  with two-reflectors of type
B are straight forward and are omitted.

It will be convenient to construct the following extensions $z^*$ 
of the function $z$ 
and $w^*$ of $w$ to the entire $\a$.
For a pair $(z,w)$ as in definition \ref{type A-B}  let
\begin{equation}\label{VQ}
V(x) = \frac{x^2}{2} + \b z(x) - \frac{\b^2}{2}~~~\mbox{and}~~~
Q(x,p) = \langle x, p \rangle + \b w(p) - \frac{p^2}{2}.
\end{equation}
It follows from  (\ref{A1}) that
\[V(x) = \sup_{p\in\Tb}Q(x,p),~~x \in \Ob .\]
That is, $V$ is convex and continuous over $\Ob$. Furthermore, since
$\Tb$ is bounded, the graph $S_V$ has no vertical supporting hyperplanes. 
For any fixed $p \in \Tb$ define the half-space
\[Q^+(p) = \{(x,Z)\in \a \times \mathbb{R}^1 ~~|~~ Z \geq Q(x,p)\}. \]
Then  
\[S_{V^*}= \partial \left ( \bigcap_{p\in\Tb}Q^+(p) \right )\]
is a graph of a convex function $V^*$ defined for all $x \in \a$. Note 
that $V^*(x) = V(x)$ when $x \in \Ob$.
We now define an extension of $z$ by putting
\[z^*(x) = \frac{1}{\b} \left [  - \frac{x^2}{2} 
+V^*(x) +  \frac{\b^2}{2} \right ]. \]

Similarly, for any fixed $x \in \Ob$ we let 
\[ R^+(x) = \{(p,Z)\in \a \times \mathbb{R}^1 ~~|~~ Z \geq R(x,p)\}, \]
where 
\[R(x,p) = \langle x, p \rangle - \b z(x) - 
\frac{x^2}{2}, ~~ p \in \a.\]
Then the function
\[U^*(p) =  \sup_{x\in\Ob}R(x,p),~~ p \in \a, \]
is defined. It is also convex.
It follows from (\ref{A2}) that for $p \in \Tb$
\[U^*(p) = \frac{p^2}{2} - \b w(p) - \frac{\b^2}{2} ~~(\equiv U(p)).\] 
The corresponding extension of $w$ we define as 
\[w^*(p) = \frac{1}{\b} \left [  \frac{p^2}{2} -U^*(p)
- \frac{\b^2}{2} \right ], ~~~p \in \a. \]
\begin{sloppypar} 
\begin{lemma} \label{lipschitz} The function $V^*$ is uniformly Lipschitz on 
$\a$ with Lipschitz
constant $\max_{\Tb}|p|$. Also, 
$U^*(p), ~p \in \a,$
is uniformly Lipschitz on $\a$ with Lipschitz
constant $\max_{\Ob}|x|$. In addition, $z \in Lip(\Ob)$ and $w \in Lip(\Tb)$
with the Lipschitz constant $ 
                \leq \frac{\sup_{(x,p) \in \Ob \times \Tb}|x-p|}{\beta}$. 
\end{lemma}
\end{sloppypar} 
\begin{proof} By our convention the normal vector to a plane $Q(x,p)$ 
(when $p$ is
fixed)
is given by $(-p,1)$. It follows from definition of $V^*$ that 
$S_{V^*}$ 
has no supporting hyperplanes with normal $(-p,1)$ such that
$p \not \in \Tb$. Since $T$ is bounded, this implies the first statement
 of the lemma.
The statements regarding $U^*$ are
established by similar arguments. From these properties of $V^*$ 
it follows that
the function $z^*$ is continuous on $\a$ and Lipschitz on any
compact subset of $\a$. Similar properties hold also for $w^*$. 

Now we estimate the Lipschitz constant for $z$ on $\Ob$. Let $(z,w)$
be a two-reflector of type A. Let $x, x'\in\Ob$ and let 
$z(x')\geq z(x)$. (If the opposite inequality holds we 
relabel $x$ and $x'$.) Fix some small $\e > 0$. It follows from 
(\ref{A1}),
(\ref{A2}) that there exists a $p'\in\Tb$ such that
$z(x') \leq k_{p',w(p')}(x') +\e$. Then
\begin{align*}
|z(x')-z(x)|& \leq k_{p',w(p')}(x')-z(x) + \e 
\leq k_{p',w(p')}(x')-k_{p',w(p')}(x) + \e\\
           & \leq \sup_{x\in\Ob}|\nabla k_{p',w(p')}(x)||x'-x| +\e
=\frac{1}{\beta}\sup_{s\in\Ob}|s-p'||x'-x| +\e\\
        &\leq \frac{1}{\beta}\sup_{s\in\Ob,p\in\Tb}|s-p||x'-x| +\e.
\end{align*}
Letting $\e \longrightarrow 0$, we obtain the statement regarding
Lipschitz constant for $z$. The same statement regarding $w$ and
two-reflectors of type B are proved similarly.
\end{proof}

Next, {\bf {we define the analogue of the ray tracing map $P$}} for
a two-reflector. For that we need to recall
the notion of the {\it {normal}} map \cite{Bak}, 
p. 114. Let $u: G \rightarrow \mathbb{R}^1$ be an arbitrary convex function defined on
domain $G \subset \a$ and $S_u$ its graph. For  $x_0 \in G$ let
\[Z-u(x_0)= \langle p,x-x_0 \rangle \]
be a hyperplane with normal $(-p,1)$ supporting to $S_u$ at $(x_0,u(x_0))$. 
The normal map $\nu_u: G \rightarrow \a$ at $x_0$ is defined as
\[\nu_u(x_0) = \bigcup\{p\},\]
where the union is taken over all hyperplanes supporting to $S_u$ at
$(x_0,u(x_0))$.

\begin{definition} Let $(z,w)$ be a two reflector of type A. 
For $x \in \a$
we put
\[
\Pt(x)=\nu_{V^*}(x).
\]
\end{definition}
For reflectors of type B the ray tracing map is defined similarly
with the use of function $-U^*$.
In general, $\Pt$ may be multivalued.  

\begin{sloppypar} 
\begin{lemma} \label{normal} Let $(z,w)$ be a two-reflector of type A
and $z^*$ and $w^*$ their respective extensions as above. 
Then $\Pt(x) \in \Tb$ for all $x \in \a$. In addition, for any $p \in \Tb$
the set $\{x \in \Ob | \Pt(x)=p\} \not = \emptyset.$
Furthermore, for any $x \in \Ob$
\begin{equation} \label{derivative}
\Pt(x) = \{\mbox{all}~~
p \in \Tb \enspace\bigl|\enspace w(p)=h_{x,z(x)}(p) \}.
\end{equation}
\begin{proof}
Let $x\in\a$ and $Q(x,p)$ a supporting hyperplane to $V^*$ at $(x,V^*(x))$.
Then the normal $p$ is in $\nu_{V^*}(x)$. On the other hand, by definition
of $V^*$, $S_{V^*}$ has only supporting hyperplanes with normals in $\Tb$. Hence,
$\Pt(x)\subset\Tb$.

Let $p\in\Tb$ and $Q(x,p)$  a supporting hyperplane to $S_{V^*}$. We need to
show that there is an $x\in\Ob$ such that $p\in\Pt(x)$. By (\ref{VQ}) and
(\ref{A2}) we have for any $x\in\Ob$
\[
	V(x)-Q(x,p)=\beta\(h_{x,z(x)}(p)-w(p) \)\geq 0.
\]
By (\ref{A2}) there exists an $x\in\Ob$ such that $V(x)-Q(x,p)=0$.
This implies the remaining two statements of the lemma.
\end{proof}
\end{lemma}
\end{sloppypar} 

\begin{remark} It follows from definition that $\Pt$ is
multivalued at points $x$ where $S_{V^*}$ has 
 more than one supporting hyperplane. At such $x$ the function $z^*$
is not differentiable.
Let $(x_0,z^*(x_0))$ be
one such point.  Then 
\[\Pt(x_0) = \{p \in \a~|~Q(x,p)~~\mbox{is supporting to $S_{V^*}$ 
at} ~(x_0,V^*(x_0))\}.\]
In other words,  a light ray
labeled by $x_0 \in \O$ that hits a point where the first reflector has a 
singular point will split into a cone of light rays. These rays will
generate a subset on the paraboloid $h_{x_0,z(x_0)}(p)$ whose projection
on $\a$ is $\Pt(x_0)$.  This is consistent with the physical 
interpretation of diffraction at
singularities of this type \cite{Keller:95}. 
\end{remark}
\begin{remark} \label{rademacher}
Since $V^*$ is
convex, by Rademacher's theorem, 
the Lebesgue measure of the set of singular points on $S_{V^*}$ is 
zero. Thus, $\Pt(x)$ is single-valued almost everywhere in $\a$.
Furthermore, the functions $z $ ($z^*$) is a difference of two convex
functions, and therefore, it is differentiable
almost everywhere in $\O$ ($\a$). The same is true for $w$ and $w^*$. 
It follows
then from the definitions of $\Pt$ and $V$  that
for almost all $x \in \O$     
\begin{equation}    \label{derivative1}
\Pt(x)= \nabla V = x+\beta \nabla z(x),
\end{equation}
A similar property holds also for the function $w$.
\end{remark}
\begin{lemma} \label{weaksol1} 
If $(z,w)$ is a two-reflector of type A 
 then for all $x\in \Ob$, $p\in \Tb$
\begin{equation} \label{equation cond on weak sol}
z(x)-w(p)\geq \frac{1}{2\beta}(\beta^2-|x-p|^2).
\end{equation}
In addition, for almost all $x\in\O$ 
there exists a unique $p\in\Tb$ such that $p = \Pt(x)$ and 
(\ref{equation cond on weak sol}) in this case is an equality.
\end{lemma}
\begin{proof}
The lemma follows from (\ref{A1}),  (\ref{A2}),
Remark~\ref{rademacher} and Lemma~\ref{normal}.
\end{proof}
Define the inverse of $\Pt$ for $p \in \Tb$
as
\[
	\Pt^{-1}(p)=\{x\in \a \enspace\bigl|\enspace p \in \Pt(x)\}.
\]
\begin{theorem}  \label{Bak1} 
Let $\cal{B}$ be the $\sigma$-algebra of Borel subsets of $T$. 
Let $(z,w)$ be a two-refllector of type A. For any set $\tau\in\cal{B}$ 
the set $\Pt^{-1}(\tau)$
is measurable relative to the standard Lebesgue measure  on $\a$. 
In addition, for any non-negative locally integrable function 
$I$ on $\a$
 the function
\[
\l(\tau)=\int_{\Pt^{-1}(\tau)} I(x)dx
\]					
is a non-negative completely additive measure on $\cal{B}$.
\end{theorem}
\begin{proof} The proof of this theorem is completely analogous
to the proofs of Theorems 9 and 16 in \cite{refl_geom1}.
\end{proof}

\begin{lemma} \label{measure preserv}
Let $\O$ and $\T$ be two bounded 
domains on $\a$ and $I$ a non-negative integrable function
on $\O$ extended to entire $\a$ by setting $I(x) \equiv 0$ for
$x \in \a \setminus \O$. 
Let $(z,w)$ be a two-reflector of type A or B.
Then for any continuous function $h$ on $\Tb$ we have 
the following change of variable formula
\begin{equation} \label{equation weak preservation of energy}
\int_{T}h(p)\l(dp)=\int_\O h(\Pt(x)) I(x)dx.
\end{equation}
\end{lemma}
\begin{proof} In the integral on the right $h(\Pt(x))$ is
discontinuous only where $\Pt$ is not single valued, that is,
on the set of measure zero. Thus, the integral on the right
is well defined.

We may assume that $\int_{\O}I(x)dx > 0$; otherwise, the statement
is trivial.
Fix some small $\e > 0$ and a positive integer N.
Partition the interval $[\min h(p), \max h(p)]$
into sub-intervals $S_1,...,S_N$ of length $ < \e/\int_{\O}I(x)dx$
and let $h_i \in S_i$. Put $\tau_i = \{p \in T~| h(p) ~ \in S_i\}$. 
Then for sufficiently large $N$
\[\mid \int_{T}h(p)\l(dp) - \sum h_i \l(\tau_i) \mid < \e.\]
For any $i,j = 1,...,N,~ i \neq j,$ 
$meas(\Pt^{-1}(\tau_i) \bigcap \Pt^{-1}(\tau_j)) =0$ 
(see Remark \ref{rademacher}).
Hence,
\[
\mid \int_{\O}h(\Pt(x))I(x)dx - \sum h_i\int_{\Pt^{-1}(\tau_i)}I(x)dx
\mid < \e. \]
This, together with the previous inequality, imply
\[
\mid \int_{\O}h(\Pt(x))I(x)dx
 - \sum h_i \l(\tau_i) \mid < 2\e.
\]

\end{proof}

\begin{definition}  A two-reflector $(z,w)$ of type A ( B)
is called a weak solution of type A (B) of 
the two-reflector problem \ref{pmain} if the map
$\Pt : \Ob \rightarrow \Tb $ is onto and
for any Borel set $\tau\subseteq T$
\[
	\l(\tau)=\int_{\tau}L(p)dp. 
\]
\end{definition}
Using Lemma \ref{measure preserv} and this definition we obtain
\begin{lemma}   \label{change of variable}
Let $(z,w)$ be a weak solution of type A  (B) 
of the two-reflector problem \ref{pmain}. Then 
for any continuous function $h$ on $\Tb$
\begin{equation}\label{weaksol}
\int_{T}h(p)L(p)dp  =\int_\O h(\Pt(x)) I(x)dx.
\end{equation}
\end{lemma}

\section{A variational problem and weak solutions of the two-reflector problem}
 \label{min_functional}

\begin{figure}
\begin{center}
\includegraphics[width=10cm]{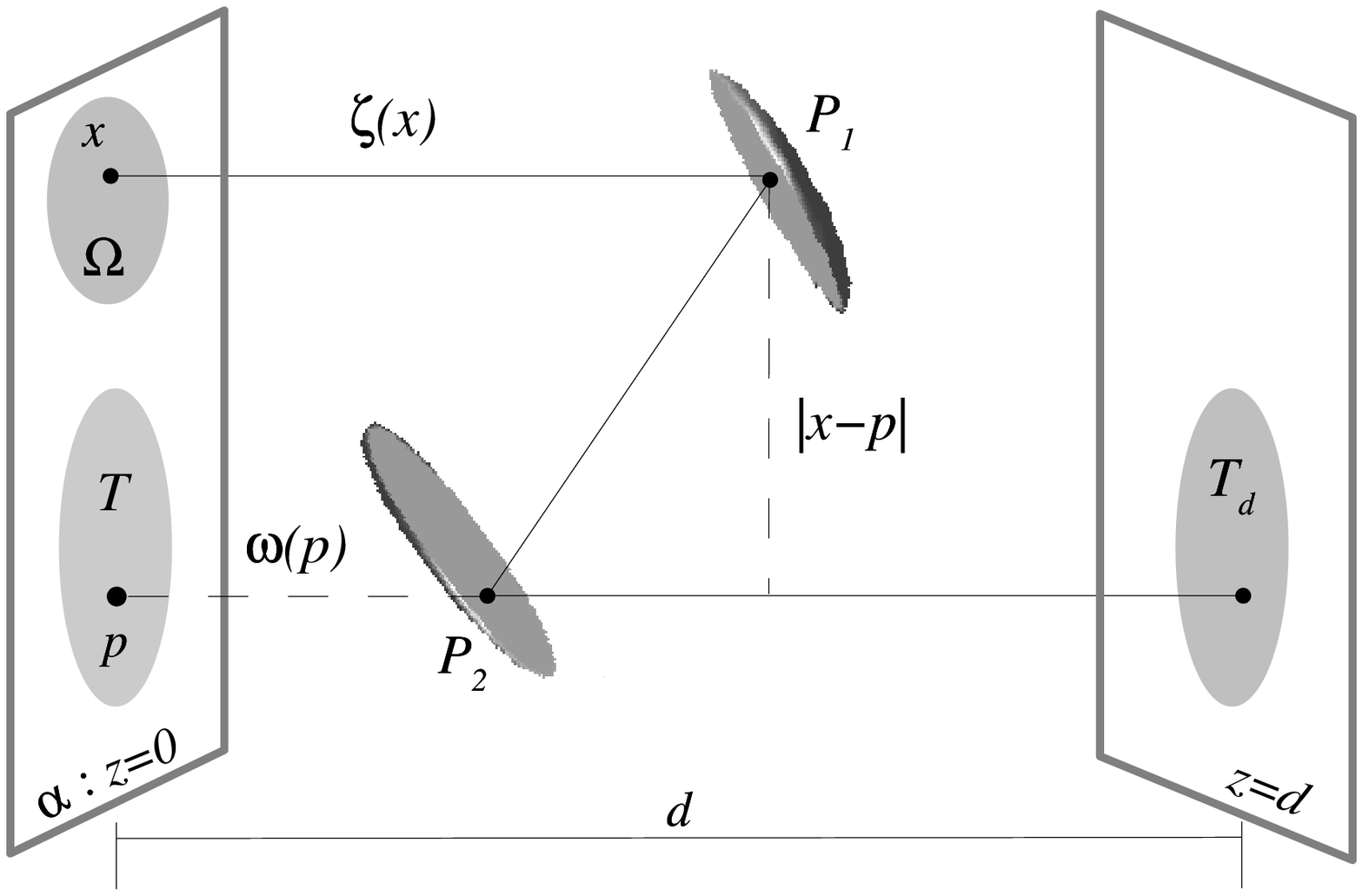}
\caption{Definition of $l\left(\z,\o,x,p\right)$.}
\label{figure bogus light}
\end{center}
\end{figure}

As before, we consider here only the case of two-reflectors of
type A. 
We comment on the case of two-reflectors of type B in section 
\ref{proof}.
Let, as before, $l$ and $d$ be the given parameters of the system, 
and $\beta=l-d$. Let 
$\zeta \in C(\Ob)$ and $\o \in C(\Tb)$. 
With any such pair $(\zeta,\o)$
we associate a ``quasi-reflector'' system in which the light
path is defined as follows. 
Let $(x,p)\in \Ob \times \Tb$ 
and let $(x,\z(x))$ be the point where the  horizontal ray
emanating from $\left(x,0\right)$
hits the graph of $\z$.  Let $P_2$ be the point where the
horizontal ray begins in order to terminate in $\left(p,d\right)$. 
 The ``light'' path is the
polygon $\overline{\left(x,0\right)P_1P_2\left(p,d\right)}$;
see Fig. \ref{figure bogus light}. The length of this path is
\begin{equation*}
l\left(\z,\o,x,p\right)=\z\left(x\right)+\sqrt{(\z(x)-
\o(p))^2+|x-p|^2}+d-\o(p).
\end{equation*}

The class of {\em admissible pairs } is defined as
\begin{equation}\label{equation def of Adm}
Adm(\O,\T)=\{ (\z,\o) \in C(\Ob) \times C(\Tb)\enspace\bigl| 
\enspace l\left(\z,\o,x,p\right) \geq l 
	\enspace \forall (x,p) \in \Ob \times \Tb\}.
\end{equation}

By construction, a pair $(\z,\o)\in C(\Ob) \times C(\Tb)$  
lies in $Adm(\O,\T)$ if and only if for all $x\in \Ob$, $p \in \Tb$
\begin{equation}\label{adm0}
\z(x)\geq k_{p,\o(p)}(x),~~\o(p)\leq h_{x,\z(x)}(p),		
\end{equation}
or, equivalently, if and only if
\begin{equation} \label{equation Adm} 
\z(x)-\o(p)\geq \frac{1}{2\beta}\left( \beta^2-|x-p|^2\right)
\end{equation}
for all $x\in \Ob$, $p \in \Tb$. It follows from (\ref{A1}), (\ref{A2}) 
that a two-reflector of type A
is an  admissible pair.

{\it {The following functional is central to our investigation.}}

Let $I$ and $L$ be two nonnegative integrable functions on $\O$ and
$T$, respectively, satisfying the energy conservation law (\ref{EC}).
For $(\z,\o)\in Adm(\O,T)$ put
\[
\F(\z,\o)=\int_\O\z(x)Idx-\int_T\o(p)Ldp.
\]

Clearly, $\F$  is linear and bounded on 
$C(\Ob)\times C(\Tb)$ with respect to
the norm $\max\{||\z||_{\infty},||\o||_{\infty}\}$.
Geometrically, $\F(\z,\o)$ is proportional to the mean horizontal distance between the points of the two graphs, the
average being weighted by the intensities.

We consider now the following 
\begin{problem} \label{problem minimize F}
	Minimize $\F$ on $Adm(\O,T)$.
\end{problem}

\begin{proposition} \cite{Brenier1, Caf_alloc:96, Gangbo/McCann:95}
\label{prop minimize F}
	There exist $(z,w)\in Adm(\O,T)$ such that 
	\[
	\F(z,w)=\inf_{(\z,\o) \in Adm(\O,T)}\F(\z,\o).
	\]
	We may further assume that the pair $(z,w)$ satisfies the conditions (\ref{A1}),
	(\ref{A2}) for a type A reflector system.

\end{proposition}
\begin{proof}
 It follows from (\ref{adm0}) that  one can restrict the search for minimizers
to two-reflectors $(\z,\o)$ of type A.
Also, because of (\ref{EC} $\F(\z,\o)$ is invariant under translations 
$\z\mapsto \z+\rho$, $\o\mapsto \o+\rho$ for
a constant $\rho\in\mathbb{R}$. Hence, it is sufficient to 
consider only two-reflectors $(\z,\o)$ for which
$\z(x_0)=0$ for some $x_0 \in \Ob$.

By Lemma \ref{lipschitz}, $\z$ and $\o$ are 
uniformly Lipschitz with the Lipschitz constant 
$K=\sup_{x\in\O,p\in\T}|x-p| / {\beta}$. 
 It follows that for all $x\in\O$
\[
|\z(x)|=|\z(x)-\z(x_0)|\leq K\diam\O.
\]
For all $p\in\Tb$ we have
\[
\o(p)\leq h_{x_0,\z(x_0)} (p)= \frac{1}{2\beta}(|x_0-p|^2-\beta^2)\leq \max_{q\in\Tb}\frac{1}{2\beta}(|x_0-q|^2-\beta^2).
\]
Finally, since $h_{x,\z(x)}(p)=(|x-p|^2-\beta^2)/2\beta+\z(x)\geq -\beta/2-K\diam\O$ for all $x\in\Ob$, we also
get the lower bound
\[
\o(p)=\inf_{x\in \O}h_{x,\z(x)}(p)\geq-\frac{\beta}{2}-K\diam\O.
\]
Therefore, the maps $\z$ and $\o$ are also uniformly bounded. 
By the Arzel\`a-Ascoli theorem, and because $\F$ is 
continuous, the infimum of $\F$ is achieved at some two-reflector
pair $(z,w)$.
\end{proof}

We now show that  weak solutions of Problem~\ref{pmain}
and solutions of the Problem~\ref{problem minimize F} are the same. 
This establishes existence of weak solutions to
Problem~\ref{pmain}.

\begin{theorem} \label{prop existence} 
Let $z \in C(\Ob),~ w \in C(\Tb)$
be a two-reflector of type A. Then
the following statements are equivalent.
	\begin{enumerate}
		\item \label{prop min} 
			$(z,w)$ minimizes $\F$ in $Adm(\O,T)$.

		\item \label{prop opt} 
			$(z,w)$ is a weak solution of type A of 
the two-reflector Problem \ref{pmain}.
	\end{enumerate}
\end{theorem}

\begin{proof}
\ref{prop opt}$\Rightarrow$\ref{prop min}.
Let $(\zeta,\omega)\in Adm(\O,T)$. Then by (\ref{equation Adm})
 and Lemma \ref{weaksol1} 
 for almost all $x\in \O$
\[
\zeta(x)-\omega(\Pt(x))\geq \frac{1}{2\beta}
\left( \beta^2-|x-\Pt(x)|^2\right)=z(x)-w(\Pt(x)).	
\]
Integrating this inequality we get
\[
\int_\O \zeta(x)Idx-\int_\O \omega(\Pt(x))Idx\geq	
\int_\O z(x)Idx-\int_\O w(\Pt(x))Idx. 
\]
Using Lemma  \ref{change of variable}, we get
\[
\F(\z,\o)=\int_\O \zeta(x)Idx-\int_T \o(p)Ldp\geq	
\int_\O z(x)Idx-\int_T w(p)Ldx=\F(z,w). 
\]
Since $(\zeta,\omega)$ was arbitrary, we are done.

\bigskip
 \ref{prop min}$\Rightarrow$\ref{prop opt}.
Let $\Pt$ denote the ray tracing map for the pair $(z,w)$. It must be shown that
	$\int_{\Pt^{-1}(\tau)}I(x)dx=\int_{\tau}L(p)dp$ for all 
Borel sets $\tau\subseteq{T}$. We will prove
	that this is the Euler-Lagrange equation for the functional $\F$.

	It is sufficient to establish this for the case when $\tau$ 
is an open ball with center $p_0 \in T$ and radius
$r>0$ contained in  $T$. For $i=1,2,\cdots$ define for $p\in\a$
\[
\chi_i(p)=\begin{cases}1,& {\text{if }} |p-p_0|<r-\frac{1}{i}\\
i(r-|p-p_0|), &{\text{if }} r-\frac{1}{i}\leq |p-p_0|<r\\
0, & {\text{if }} |p-p_0|\geq r.
\end{cases}
\]
Then $\chi_i$ is continuous on $\a$, $0\leq\chi_i\leq1$, and the sequence 
$\{\chi_i(p):p\in\a\}_{i=1}^{\infty}$ converges on $\a$ pointwise
	to the characteristic
        function of $\tau$  $\chi_\tau(p)$ .
	
	Fix some $i$ and for  $\e \in (-1,1)$ put 
	\begin{align*}
	&w_{\e}(p)=w(p)+\e\cdot\chi_i(p)\\
	&z_{\e}(x)=\sup_{p\in\Tb}k_{p,w_\e(p)}(x)
			=\sup_{p\in\Tb}\{\frac{1}{2\beta}\left( \beta_{0}^2-|x-p|^2\right)+w_\e(p)\}.			
	\end{align*}
By construction, the pair $(z_\e,w_\e)$ satisfies 
the condition (\ref{equation Adm}), with
$\z$ replaced by $z_\e$ and $\o$ repaced by $w_\e$. 
We show now that $z_\e$ belongs to $\Lip(\Ob)$.
Let $x, x'\in\Ob$ and
let $z(x')\geq z(x)$. (If the opposite inequality holds we relabel $x$ and $x'$.)  
Let $p'\in\Tb$ be such that
$z_{\e}(x')=k_{p',w_\e(p')}(x')$. Then
	\begin{align*}
	|z_{\e}(x')-z_{\e}(x)|&=k_{p',w_{\e}(p')}(x')-z_{\e}(x)\\
		&\leq k_{p',w_{\e}(p')}(x')-k_{p',w_{\e}(p')}(x) \\
	   & \leq \sup_{s\in\Ob}|\nabla k_{p',w_{\e}(p')}(s)||x'-x|=\frac{1}{\beta}\sup_{s\in\Ob}|s-p'||x'-x| \\
		&\leq \frac{1}{\beta}\sup_{s\in\Ob,p\in\Tb}|s-p||x'-x|.
	\end{align*}  
Hence, $z_\e$ is continuous  and $(z_\e,w_\e)\in Adm(\O,\T)$.

	Now let $x\in\Ob$. For each $\e$ let $p_\e$ be a point in $\Tb$ such that
$z_{\e}(x)=k_{p_{\e},w_{\e}(p_{\e})}(x)$. (This choice, of course,
may not be unique.) Then
	\begin{align*}
		z_\e(x)-z(x)&=k_{p_{\e},w_{\e}(p_{\e})}(x)-z(x)\leq k_{p_{\e},w_{\e}(p_{\e})}(x)-k_{p_{\e},w(p_{\e})}(x)\\
			&=w_{\e}(p_{\e})-w(p_\e)=\e \chi_i(p_\e).
	\end{align*}
	Similarly, if $p\in\Pt(x)$, then
	\begin{align*}
		z_\e(x)-z(x)&=z_{\e}(x)-k_{p,w(p)}(x)\geq  k_{p,w_{\e}(p)}(x)-k_{p,w(p)}(x)\\
				&=\e\chi_i(p).
	\end{align*}
Therefore,
\begin{equation}    \label{ineq for z}
-|\e|\leq \e\chi_i(p)\leq z_\e(x)-z(x)\leq \e\chi_i(p_\e)\leq |\e|
	\end{equation}
for all $x\in\Ob$. In particular, $z_\e$ converges uniformly to $z$ on $\Ob$ as $\e\rightarrow 0$.
	
Now consider those $x\in\O$ for which the ray tracing
map $\Pt$ is single-valued. This is the case for almost all 
$x\in\O$. We claim that $p_{\e} =\Pt_{\e}$, where
$\Pt_{\e}$ denotes the ray tracing map for $(z_{\e}, w_{\e})$,
 converge  to $p$ as $\e\rightarrow 0$.

	Suppose this is not true. Then there is a sequence $\{p_{\e_{j}}\}$, $j=1,2,\dots$ with 
	$\e_{j}\rightarrow 0$ as $j\rightarrow\infty$ and a constant
	$\eta>0$ such that
	\[
		|p-p_{\e_{j}}|>\eta\text{ for all }j.
	\]
	Let $z'(x)=\max_{p'\in\Tb, |p'-p|\geq \eta}k_{p',w(p')}(x)$. 
Note that the maximum in the definition of
	$z'$ is attained. Since $p$ is the unique point in $\Tb$ such that $z(x)=k_{p,w(p)}(x)$, it follows
	that $z'<z(x)$. Therefore, for all $j$
	\begin{align*}
		z(x)-z_{\e_j}(x)&=z(x)-k_{p_{\e_{j}},w_{\e_{j}}(p_{\e_{j}})}(x)
			=z(x)-k_{p_{\e_{j}},w(p_{\e_{j}})}(x)+w(p_{\e_{j}})-w_{\e_{j}}(p_{\e_{j}})\\
			&=z(x)-k_{p_{\e_{j}},w(p_{\e_{j}})}(x)-\e_{j}\chi_i(p_{\e_{j}})\\
			&\geq z(x)-z'-|\e_j|.	
	\end{align*}
	 This contradicts to the fact that $z_{\e_{j}}(x)$ converges to $z(x)$ as 
	$j\rightarrow\infty$. 
Hence, $p_\e\rightarrow p$ if $\e\rightarrow 0$.

It follows from (\ref{ineq for z}) that
\[
\Bigl|\frac{z_\e(x)-z(x)}{\e}-
\chi_i(p)\Bigr|\leq |\chi_i(p_\e)-\chi_i(p)|.
\]	
Letting $\e\rightarrow 0$ and using the continuity of $\chi_i$,
we conclude that for almost all $x \in \O$
	\[
			\frac{d}{d\e}\bigg|_{\e=0}z_\e(x)=\chi_i(p)=\chi_i(\Pt(x)).
	\]
Then
	\[
			\frac{d}{d\e}\bigg|_{\e=0}\int_{\O} z_\e(x)I(x)dx=\int_{\O} \chi_i(\Pt(x)) I(x)dx.
	\]
Since $\F$ has a minimum at $(z,w)$, we obtain
	\[
		0=\frac{d}{d\e}\bigg|_{\e=0}\F(z_\e,w_\e)=\int_{\O}\chi_i(\Pt(x))I(x)dx-\int_{\T} \chi_i(p)L(p)dp.
	\]
	Now let $i\rightarrow\infty$ in this equality. 
This is possible as $\chi_i(p)\rightarrow\chi_\tau(p)$ pointwise
on $\a$ and therefore $\chi_i(\Pt(x))\rightarrow\chi_\tau(\Pt(x))$ 
pointwise almost everywhere on $\O$. Then, noting that
for almost all $ x \in \O~~
 \chi_{\tau}(\Pt(x)) = \chi_{\Pt^{-1}(\tau)}(x)$, we obtain
	\begin{align*}
		\int_{\tau}L(p)dp&=
\int_{\T}\chi_{\tau}(p)L(p)dp=\int_{\O}\chi_{\tau}(\Pt(x))I(x)dx\\
				&=\int_{\O}\chi_{\Pt^{-1}(\tau)}(x)I(x)dx=\int_{\Pt^{-1}(\tau)}I(x)dx.		
	\end{align*}
	This completes the proof of the theorem.
\end{proof}

\section{Connection between the
two-reflector and Monge-Kantorovich problems} \label{MK}

In our notation, the Monge-Kantorovich mass transfer problem 
\cite{Brenier1}
can be formulated as follows. Consider the class of maps 
$P:\O\rightarrow T$ which are measure-preserving, that is, 
they satisfy the substitution rule
\[
\int_{\O}h(P(x))Idx=\int_Th(p)Ldp
\]
for all continuous functions $h$ on $\Tb$. Each such map is
called a {\em plan}.
 

\begin{problem} Minimize the quadratic transportation cost
\begin{equation} \label{equ transp cost}
P\mapsto \frac{1}{2}\int_{\O}|x-P(x)|^2Idx.
\end{equation}
among all  plans $P$.
\end{problem}

Note that for weak solutions to the  two-reflector 
Problem~\ref{pmain} the ray tracing map $\Pt$ 
is a plan by Lemma~\ref{change of variable}. 
In fact, we have the following
\begin{theorem} \label{prop mk} 
Let $(z,w)$ be a weak solution of type A of the two-reflector 
Problem \ref{pmain}.
	Let $\Pt$ be the corresponding ray tracing map.
	Then $\Pt$ minimizes the quadratic transportation cost 
	(\ref{equ transp cost}) among all plans $P:\O\rightarrow T$, and any other minimizer is equal to $\Pt$
	almost everywhere on $\supp(I)\setminus \{x\in\Ob:I(x)=0\}$.
\end{theorem}
\begin{proof}
	Let $P:\O\rightarrow T$ be any plan. Then by
	(\ref{equation cond on weak sol}),
	\begin{equation}       \label{eqn plan inequality}
	z(x)-w(P(x))\geq \frac{1}{2\beta}(\beta^2-|x-P(x)|^2),
	\end{equation}
	for almost all $x\in \O$ and equality holds if and only if $P(x)=\Pt(x)$.
	Integrating against $Idx$  and applying Lemma 
\ref{change of variable}, we get
	\begin{align*} 
	\frac{1}{2 \beta}\int_{\O}[\beta^2-&|x-P(x)|^2]Idx \leq\int_{\O}[ z(x)-w(P(x))]Idx      \\
	&=\int_{\O}z(x)Idx-\int_{\T}w(p)Ldp 
	=\int_{\O}[ z(x)-w(\Pt(x))]Idx  \\
	&=\frac{1}{2\beta}\int_{\O}[\beta^2-|x-\Pt(x)|^2]Idx.
	\end{align*}
This shows that $\Pt$ is a minimizer of the transportation cost. 

To show uniqueness, note that if equality holds
in the integral inequality, then equality must hold in 
(\ref{eqn plan inequality}) 
for almost all $x\in\supp{(I)}\setminus \{I=0\}$. 
Therefore, $P\equiv \Pt$ a.e. on $\supp(I)\setminus \{I=0\}$.
\end{proof} 

\begin{remark}
The functional $\F$ is the dual of the quadratic cost functional 
(\ref{equ transp cost}) \cite{Brenier1}. 
\end{remark}

\section{Existence and uniqueness of weak solutions to the two-reflector problem}
 \label{proof}
\begin{theorem}	\label{thm existence & uniqueness}
There exist weak solutions of type A and of type B to 
the two-reflector Problem \ref{pmain}.
Furthermore, if $(z,w)$ and $(z',w')$ are two solutions of the 
same type with ray tracing maps $\Pt$ and $\Pt'$, respectively, then
\[
		\Pt(x)\equiv\Pt'(x)
\]
for almost all $x\in\supp(I)\setminus \{x\in\O:I(x)=0\}$.
\end{theorem}
\begin{proof} By Proposition \ref{prop minimize F} and 
Theorem \ref{prop existence}
we know that Problem \ref{pmain} has a solution. The only property that
remains to be checked is that for that solution, $\Pt\colon\Ob\rightarrow\Tb$
is onto. But this follows from Lemma~{\ref{normal}}.

In order to prove uniqueness, 
let $(z,w)$ and $(z',w')$ be two solutions of type A, with
ray tracing maps $\Pt$ and $\Pt'$, respectively. Then by Proposition
\ref{prop mk}, both are minimizers of the quadratic Monge-Kantorovich
cost functional, so that $\Pt=\Pt'$ a.e. on $\supp(I)\setminus
\{I=0\}$.

So far, existence and uniqueness of weak solutions to
the two-reflector Problem \ref{pmain} has been shown for
weak solutions
of type A.
However, it is clear that a similar result holds for weak
solutions of type B as
well. Namely, we can define $Adm_+(\O,T)$ as the space of all pairs
of reflectors such that $l(\xi,\omega,x,p)$ is {\em less} than $l$,
and then that $\F$ admits a {\em maximizing} pair on $Adm_+(\O,T)$,
which is a weak solution of the two-reflector problem~\ref{pmain}.
This shows, in particular, that for such a solution the ray tracing map
maximizes the quadratic transportation cost among all plans.
\end{proof}

\begin{corollary}
	Suppose that in addition to the assumptions in 
Theorem \ref{thm existence & uniqueness} the function
	$I > 0 $ in $\O$.
	Then there is a constant
	$\rho\in {\mathbb R}$ such that
	\begin{align*}
	& z'(x)=z(x)+\rho\\
	& w'(p)=w(p)+\rho		
	\end{align*}
	for all $x\in\Ob$ and all $p\in\Tb$. In other words, weak solutions of type A
	 are unique on $\Ob$ and $\Tb$ up
	 to a
	translation of the reflector system, and the same result holds for type B solutions.
\end{corollary}
\begin{proof} We show this for weak solutions of type A. By the theorem,
$\Pt(x)\equiv\Pt'(x)$ for {\em all}
$x\in\O$. Note that by (\ref{derivative1}),  $\nabla z(x)=\nabla z'(x)$ for almost all
$x\in\O$. It follows that there is a constant $\rho$ such that $z'(x)=z(x)+\rho$ for all $x\in\Ob$. 
Now, by
definition of $w'(p)$,
\[
w'(p)=\inf_{x\in\Ob}h_{x,z'(x)}(p)=\inf_{x\in\Ob}h_{x,z(x)}(p)+\rho=w(p)+\rho
\]
for all $p\in\Tb$.
\end{proof}

\bibliographystyle{plain}
\bibliography{main}

\end{document}